# On the Structure of Some One-generator Braces

L.A. Kurdachenko, I.Ya. Subbotin

**Abstract**. We describe the one-generator braces A satisfying the condition $A^3 = <0>$.

**Mathematics Subject Classification**: Primary: 17D99: None of the above but in this section.

Secondary:16T25: Yang-Baxter equations, 81R50: Quantum groups and related algebraic methods.

Acknowledgments: The first author is grateful for the support Isaac Newton Institute for Mathematical Sciences and the University of Edinburgh provided in the frame of LMS Solidarity Supplementary Grant Program.
The authors are sincerely grateful to Professor A. Smoktunowicz for useful consultations.

A *left brace* is a set A with two binary operations + and · satisfying the following conditions: A is an abelian group by addition, A is a group by multiplication, and $a(b + c) = ab + ac - a$ for every $a, b, c \in A$.

In order to help study involutive set-theoretic solutions of the Yang-Baxter equation, W. Rump [RW2005, RW2007A] introduced braces as a generalization of Jacobson radical rings. The main reason for introducing this algebraic object is it allowing another possible strategy to attack the problem of classifying such solutions. Subsequent papers [RW2007B, RW2007C, RW2008, RW2009] found connections of braces with other algebraic structures. Currently, the theory of braces is developing very intensely. An important part of this theory is the study of the internal structure of braces. One of the immediately arising questions here is the study of braces generated by one element (i.e., one-generator braces). In general, this task is very difficult, so the study of a one-generator brace should be started under some additional conditions:

A *left brace* is a set A with two binary operations + and · satisfying the following conditions: A is an abelian group by the addition, A is a group by multiplication, and $a(b + c) = ab + ac - a$ for every $a, b, c \in A$.

Let A be a left brace. A subset S of A is called *subbrace* (more precisely *left subbrace*) if S is closed by the addition and multiplication and is the left braces by its restriction on S (i.e., S by the restriction of addition is a subgroup of an additive group of A and by the restriction of multiplication, is a subgroup of a multiplicative group of A).

Let A be a left brace. For every element $a \in A$, we define the function $\lambda_a$: $A \longrightarrow A$ by the rule $\lambda_a(x) = ax - a$ for all elements $x \in A$.

Let A be a left brace. Put $a*b = ab - a - b$. We can see that $a*b = \lambda_a(b) - b$. This new operation plays a very important role in left braces.

A left brace A is called *trivial or abelian* if $a*b = 0$ or $a + b = ab$ for all elements $a, b \in A$.

The subbrace L of a brace A is called an *ideal* if $a * z, z * a \in L$ for all elements $a \in A$ and $z \in L$.

If K, L are subbraces of A, then denote by $K * L$ the subgroup of the additive group of A generated by the elements $x*y$, where $x \in K, y \in L$.

Let A be the left brace. Put $A^{(1)} = A$ and, recursively, $A^{(\alpha+1)} = A^{(\alpha)} * A$ for all of ordinal $\alpha$ and $A^{(\lambda)} = \bigcap_{\mu < \lambda} A^{(\mu)}$ for limit ordinal $\lambda$. We note that $A^{(\alpha)}$ is an ideal of A for every ordinal $\alpha$.

A subbrace L of a left brace A is called *a left ideal of A* if $a * b \in L$ for every element $a \in A$ and every element $b \in L$.

Put $A^1 = A$ and, recursively, $A^{\alpha+1} = A * A^{(\alpha)}$ for all of ordinal $\alpha$ and $A^\lambda = \bigcap_{\mu < \lambda} A^\mu$ for limit ordinals $\lambda$. We note that $A^\alpha$ is a left ideal of A for every ordinal $\alpha$.

Let A be a left brace and let $\mathfrak{S}$ be a family of subbraces of A. Then the intersection S of all subbraces of family $\mathfrak{S}$ is the subbrace of A. If M is the subset of A, then let $\mathfrak{M}$ be the family of all subbraces of A including M. Then the intersection of all subbraces of family $\mathfrak{M}$ is the least subbrace of A including M. This subbrace is called the *subbrace of A generated by a subset M* and will be denoted by $\mathbf{br}(M)$.

A subbrace S of a brace A is called *finitely generated* if there exists a finite subset M such that $A = \mathbf{br}(M)$. If $M = \{a\}$, then A is called a *one-generator brace.*

If $M = \{a\}$, then we will talk about the *one-generator subbrace of A.*

One of the first basic steps pertaining to the study of left braces is the study of the one-generator braces. In general, it is a complicated problem. Every case of such study is important. In this paper, we will study one-generator braces A such that $A^3 = <0>$.

Like in other algebraic structures, there is the concept of nilpotency in braces. In the theory of braces, there are different approaches to this concept (see, for example, papers [**CF2018**, **CSV2019**, **JvAV2022**, **SA2018**, **CGIS2017**]). In the case of a one generator brace such that $A^3 = <0>$, we come to the following notion of nilpotency:

We say that A is a series that is called *nilpotent in the sense of Smoktunowicz* if there are positive integers n, k such that $A^{(n)} = <0> = A^k$. These braces have been introduced in a paper of A. Smoktunowicz [**SA2018**]. Denote by $\mathcal{N}_S(n, k)$ the class of left braces satisfying $A^{(n)} = <0> = A^k$ where n, k are the least integers that contain this property.

The main results of the current paper are the following:

**THEOREM A.** *Let A be a left brace such that $A^3 = <0>$. Suppose that A is generated by the element a, and put $a_1 = a, a_2 = a * a = a_1 * a, a_3 = a_2 * a$. Then every element of A has a form*

$k_1a_1 + k_2a_2 + k_3a_3$, where $k_1, k_2, k_3$ are integers. Moreover, if $x = k_1a_1 + k_2a_2 + k_3a_3$, $y = t_1a_1 + t_2a_2 + t_3a_3$ are the elements of $A$, $k_1, k_2, k_3, t_1, t_2, t_3 \in \mathbb{Z}$, then

$$x \star y = t_1k_1a_2 + \tfrac{1}{2}(2k_2 + k_1 - k_1^2)t_1a_3 \text{ and}$$
$$xy = (k_1 + t_1)a_1 + (k_2 + t_2 + t_1k_1)a_2 + (k_3 + t_3 + \tfrac{1}{2}(2k_2 + k_1 - k_1^2)t_1)a_3.$$

In particular, if $a_3 = 0$, then every element of $A$ has a form $k_1a_1 + k_2a_2$, $k_1, k_2$ are integers. Moreover, if $x = k_1a_1 + k_2a_2$, $y = t_1a_1 + t_2a_2$ are elements of $A$, $k_1, k_2, t_1, t_2 \in \mathbb{Z}$, then

$$x \star y = t_1k_1a_2 \text{ and } xy = (k_1 + t_1)a_1 + (k_2 + t_2 + t_1k_1)a_2.$$

Furthermore, $A$ is nilpotent in the sense of Smoktunowizc. More precisely, $A \in \mathcal{N}_S(3, 3)$ whenever $a_3 = 0$, and $A \in \mathcal{N}_S(4, 3)$ whenever $a_3 \neq 0$.

**THEOREM B1.** *There exists a one-generator left brace $D = \mathbb{Z} \times \mathbb{Z}$ where the addition and multiplication are defined by the rules:*

$$(k_1, k_2) + (t_1, t_2) = (k_1 + t_1, k_2 + t_2),$$
$$(k_1, k_2)(t_1, t_2) = (k_1 + t_1, k_1t_1 + k_2 + t_2).$$

*If $A$ is an arbitrary one-generator left brace such that $A \in \mathcal{N}_S(3, 3)$, then there exists an epimorphism $f: D \longrightarrow A$.*

**THEOREM B2.** *There exists a one-generator left brace $D = \mathbb{Z} \times \mathbb{Z} \times \mathbb{Z}$ where the addition and multiplication are defined by the rules:*

$$(k_1, k_2, k_3) + (t_1, t_2, t_3) = (k_1 + t_1, k_2 + t_2, k_3 + t_3),$$
$$(k_1, k_2, k_3)(t_1, t_2, t_3) = (k_1 + t_1, k_1t_1 + k_2 + t_2, k_3 + t_3 + \tfrac{1}{2}(2k_2 + k_1 - k_1^2)t_1).$$

*If $A$ is an arbitrary one-generator left brace such that $A \in \mathcal{N}_S(4, 3)$, then there exists an epimorphism $f: D \longrightarrow A$.*

Thus, we can see that the left brace, constructed in **Theorem B1** (respectively, in **Theorem B2**), is a free one-generator left brace such that $A^3 = <0> = A^{(3)}$ (respectively, $A^3 = <0> = A^{(4)}$).

## 1. Some preliminary results

We will also need the following properties of operation $\star$ and the mapping $\lambda_a$. One can find the proofs of these results in papers [**CF2018**, **JvAV2022**], for example.

**1.1. LEMMA.** *Let $A$ be a left brace. Then*

$$a \star (b + c) = a \star b + a \star c, \ (ab) \star c = a \star (b \star c) + b \star c + a \star c,$$
$$(a + b) \star c = a \star (\lambda_{a^{-1}}(b) \star c) + (\lambda_{a^{-1}}(b) \star c) + a \star c.$$
$$\lambda_y(b \star a) = yby^{-1} \star \lambda_y(a), \ yby^{-1} = \lambda_y(\lambda_b(y^{-1}) - y^{-1} + b) = \lambda_y(b \star y^{-1} + b)$$

*for all elements $a, b, c, y \in A$.*

The following result, the proof of which is possible to find in paper [**CF2018**], will also be used by us.

**1.2. PROPOSITION.** *Let A be a left brace and L be a left ideal of A. Then $L * A$ and $A * L$ are left ideals of A. Moreover, if L is an ideal of A, then $L * A$ is an ideal of A.*

Using **Proposition 1.2,** we obtain:

**1.3. PROPOSITION.** *Let A be a left brace. Then $A^\alpha$ is a left ideal for each ordinal $\alpha$, and $A^{(\alpha)}$ is an ideal for each ordinal $\alpha$.*

A brace A is called *soluble* if A has a finite series

$$<0> = A_0 \leq A_1 \leq \ldots \leq A_{n-1} \leq A_n = A$$

of ideals whose factors $A_j/A_{j-1}$ are abelian, $1 \leq j \leq n$.

**1.4. PROPOSITION.** *Let A be a left brace and suppose that $A^3 = <0>$. Then*

*$A^2$ is abelian. Moreover, $x * y = 0$ ( or $xy = x + y$ ) for every element $x \in A$, $y \in A^2$,*
*$(xz) * y = x * y + z * y$ for every element $x, y, z \in A$,*
*$g^n * z = n(g * z)$ for every integer n and every element $g, z \in A$.*

**PROOF.** The equality $A^3 = <0>$ implies that $x * y = 0$ for every element $x \in A$, $y \in A^2$. It follows that $xy = x + y$. In particular, it follows that $A^2 = A^{(2)}$ is abelian.

By **Lemma 1.1,** for every element $x, y, z$ of A, we have:

$$(xz) * y = x * (z * y) + x * y + z * y = x * y + z * y.$$

Furthermore,

$$0 = 0 * z = 1 * z = (g^{-1} g) * z = g * z + g^{-1} * z.$$

Hence, $g * z + g^{-1} * z = 0$. It follows that $g^{-1} * z = - g * z$. Furthermore,

$$g^2 * z = (gg) * z = g * z + g * z = 2(g * z).$$

Using ordinary induction, we obtain that $g^n * z = n(g * z)$ for every positive integer n.

Now, let n be a negative integer. Then $n = - k$, $k > 0$. We have $g^n = (g^{-1})^k$. By what is proven above,

$$g^n * z = ((g^{-1})^k * z) = k((g^{-1}) * z) = -k(g * z) = n(g * z).$$

**1.5. PROPOSITION.** *Let A be a left brace and suppose that $A^3 = <0>$. Then*

$$a^m = ma + \tfrac{1}{2}(m^2 - m)(a * a)$$

*for every element a of A and every integer m.*

**PROOF.** Put $b = a * a$. From $a * a = a^2 - a - a = a^2 - 2a$, we obtain that $a^2 = 2a + b$. Further,

$$a^3 = aa^2 = a(a + a + b) = a^2 + a^2 + ab - 2a = 2a + b + 2a + b + ab - 2a = 2a + 2b + ab.$$

By what is noted above, $ab = a + b$, so that $a^3 = 2a + 2b + a + b = 3a + 3b$. Again, we have:

$$a^4 = aa^3 = a(a + a + a + b + b + b) = 3a^2 + 3ab - 5a = 3(2a + b) + 3a + 3b - 5a = 4a + 6b,$$

$$a^5 = aa^4 = a(4a + 6b) = 4a^2 + 6ab - 9a = 4(2a + b) + 6a + 6b - 9a = 5a + 10b,$$
$$a^6 = aa^5 = a(5a + 10b) = 5a^2 + 10ab - 14a = 5(2a + b) + 10a + 10b - 14a = 6a + 15b.$$

Applying ordinary induction, we obtain that $a^m = ma + ½(m^2 - m)b$ for arbitrary positive integer m.

Now, suppose that m < 0. Let $y = a^{-1}$. Then $a^m = y^t$, where $t = -m > 0$. By what is proven above, we have $a^m = y^t = ty + ½ (t^2 - t)z_1$, where $z_1 = y * y = (a^{-1}) * (a^{-1})$. Using what is stated above, we obtain $(a^{-1}) * (a^{-1}) = - (a) * (a^{-1})$. Now, we have:

$$a * (a^{-1}) = aa^{-1} - a - a^{-1} = 1 - a - a^{-1} = 0 - a - a^{-1} = - a - a^{-1},$$
$$(a^{-1}) * a = a^{-1}a - a - a^{-1} = - a - a^{-1},$$

so that $a * (a^{-1}) = (a^{-1}) * a = - a * a = - b$. Hence, $(a^{-1}) * (a^{-1}) = a * a = b$ and

$$a^m = y^t = ty + ½ (t^2 - t)b = t(a^{-1}) + ½ (t^2 - t)b.$$

As we have seen above, $(a^{-1}) * a = - a - a^{-1}$. On the other hand, $(a^{-1}) * a = -(a * a) = -b$, so that $- a - a^{-1} = -b$ and $a^{-1} = b - a$. Thus, $t(a^{-1}) = t(b - a)$ and

$$a^m = y^t = t(a^{-1}) + ½ (t^2 - t)b = t(b - a) + ½ (t^2 - t)b = t(-a) + ½ (t^2 + t)b =$$
$$(-t)a + ½ (t^2 + t)b = ma + ½ (m^2 - m)b.$$

**1.6. PROPOSITION.** *Let A be a left brace and suppose that $A^3 = < 0 >$. Then*

$$(na) * (ka) = kn(a * a) + ½ (n - n^2)k((a * a) * a)$$

*for every element a of A and every integer n, k.*

**PROOF.** Put $b = a * a$, $c = b * a = (a * a) * a$. We have:

$$a^n + a^{-n} = na + ½ (n^2 - n)b + (-n)a + ½ (n^2 + n)b = n^2b.$$

It follows that $a^n + a^{-n} - n^2b = 0$, and therefore $- a^n = a^{-n} - n^2b$.

Now, We will find the element $(na) * (ka)$. Using **Lemma 1.1,** we obtain that $(na) * (ka) = k((na) * a)$. Equality $a^n = na + ½ (n^2 - n)b$ implies that

$$na = a^n - ½ (n^2 - n)b = a^n + ½ (n - n^2)b.$$

Using **Proposition 1.5,** we obtain that

$$a^n + ½ (n - n^2)b = a^n (½ (n - n^2)b).$$

An application of **Lemma 1.1** and **Proposition 1.5** give

$$(na) * a = (a^n (½ (n - n^2)b)) * a = a^n * a + (½ (n - n^2)b) * a = n(a * a) + ½ (n - n^2)(b * a).$$

Thus,

$$(na) * (ka) = knb + ½ (n - n^2)kc.$$

## 2. The structure of a one-generator brace A such that $A^3 = < 0 >$

**2.1. PROPOSITION.** *Let A be a left brace such that $A^3 = <0>$. Let a be an element of A and put $a_1 = a$, $a_2 = a * a = a_1 * a$, $a_3 = a_2 * a$, $a_{j+1} = a_j * a$, $j \in \mathbf{N}$. Then the subbrace of A, generated by element a, is a subset of elements having the following form: $\sum_{j \in \mathbf{N}} k_j a_j$, $k_j \in \mathbf{Z}$ and $k_j \neq 0$ only for finitely many indexes $j \in \mathbf{N}$.*

**PROOF.** We note that the additive subgroup of A generated by the subset $\{a_j \mid j \in \mathbf{N}\}$ consists of the elements $\sum_{j \in \mathbf{N}} k_j a_j$ where $k_j \neq 0$ only for finitely many indexes $j$.

Let $x = \sum_{j \in \mathbf{N}} k_j a_j$, $y = \sum_{j \in \mathbf{N}} t_j a_j$. Using **Lemma 1.1** and **Proposition 1.4**, we obtain

$$x * y = x * (\sum_{j \in \mathbf{N}} t_j a_j) = \sum_{j \in \mathbf{N}} (x * (t_j a_j)) = \sum_{j \in \mathbf{N}} t_j ((x * a_j) = t_1 (x * a_1).$$

**Proposition 1.4** implies that $\sum_{j \in \mathbf{N}} k_j a_j = (k_1 a_1)(\sum_{j > 1} k_j a_j)$, and the fact that $A^2$ is abelian implies that

$$k_j a_j = a_j^{k_j}, j > 1, \sum_{j > 1} k_j a_j = \prod_{j > 1} k_j a_j = \prod_{j > 1} a_j^{k_j},$$

and we obtain

$$x * a_1 = (\sum_{j \in \mathbf{N}} k_j a_j) * a_1 = (k_1 a_1)(\sum_{j > 1} k_j a_j) * a_1 = (k_1 a_1) * a_1 + (\sum_{j > 1} k_j a_j) * a_1 =$$
$$(k_1 a_1) * a_1 + (\prod_{j > 1} a_j^{k_j}) * a_1 = (k_1 a_1) * a_1 + \sum_{j > 1} (a_j^{k_j} * a_1) =$$
$$(k_1 a_1) * a_1 + \sum_{j > 1} k_j (a_j * a_1) = (k_1 a_1) * a_1 + \sum_{j > 1} k_j a_{j+1}.$$

By **Proposition 1.6**, $(k_1 a_1) * a_1 = k_1 a_2 + \frac{1}{2}(k_1 - k_1^2) a_3$, so that

$$x * a_1 = k_1 a_2 + \frac{1}{2}(k_1 - k_1^2) a_3 + k_2 a_3 + k_3 a_4 + \sum_{j > 3} k_j a_{j+1} =$$
$$k_1 a_2 + \frac{1}{2}(2k_2 + k_1 - k_1^2) a_3 + k_3 a_4 + \sum_{j > 3} k_j a_{j+1}.$$

And we have:
$$x * y = t_1 (x * a_1) = t_1 k_1 a_2 + \frac{1}{2}(2k_2 + k_1 - k_1^2) t_1 a_3 + k_3 t_1 a_4 + \sum_{j > 3} k_j t_1 a_{j+1}.$$

By $x * y = xy - x - y$ we obtain

$$xy = (\sum_{j \in \mathbf{N}} k_j a_j)(\sum_{j \in \mathbf{N}} t_j a_j) =$$
$$t_1 k_1 a_2 + \frac{1}{2}(2k_2 + k_1 - k_1^2) t_1 a_3 + k_3 t_1 a_4 + \sum_{j > 3} k_j t_1 a_{j+1} + \sum_{j \in \mathbf{N}} k_j a_j + \sum_{j \in \mathbf{N}} t_j a_j =$$
$$(k_1 + t_1) a_1 + (k_2 + t_2 + t_1 k_1) a_2 + (k_3 + t_3 + \frac{1}{2}(2k_2 + k_1 - k_1^2) t_1) a_3 + \sum_{j \geq 4} (k_j + t_j + t_1 k_{j-1}) a_j.$$

It is possible to check that

$$x^{-1} = -k_1 a_1 + (k_1^2 - k_2) a_2 + (\frac{1}{2}(2k_2 + k_1 - k_1^2) k_1 - k_3) a_3 + t_3 + \frac{1}{2}(2k_2 + k_1 - k_1^2) t_1) a_3 +$$
$$\sum_{j \geq 4} (k_1 k_{j-1} - k_j) a_j.$$

Thus, we can see that the set of elements having a form $\sum_{j \in \mathbf{N}} k_j a_j$, $k_j \in \mathbf{Z}$ for all indexes $j \in \mathbf{N}$ is a subbrace. On the other hand, every subbrace containing an element $a$ contains the elements $\sum_{j \in \mathbf{N}} k_j a_j$, $k_j \in \mathbf{Z}$ for all indexes $j \in \mathbf{N}$. Hence, we obtain that a subbrace, generated by element $a$ coincides with the set of elements having a form $\sum_{j \in \mathbf{N}} k_j a_j$, $k_j \in \mathbf{Z}$ for all indexes $j \in \mathbf{N}$.

### PROOF OF THEOREM A.

By **Proposition 2.1**, every element $x$ of A has a form $x = \sum_{j \in \mathbf{N}} k_j a_j$, $k_j \in \mathbf{Z}$ for all indexes $j \in \mathbf{N}$. From the proof of **Proposition 2.1**, we also obtain that

$$(\sum_{j \in \mathbf{N}} k_j a_j)(\sum_{j \in \mathbf{N}} t_j a_j) =$$
$$t_1 k_1 a_2 + \tfrac{1}{2}(2k_2 + k_1 - k_1^2)t_1 a_3 + k_3 t_1 a_4 + \sum_{j > 3} k_j t_1 a_{j+1} + \sum_{j \in \mathbf{N}} k_j a_j + \sum_{j \in \mathbf{N}} t_j a_j =$$
$$(k_1 + t_1)a_1 + (k_2 + t_2 + t_1 k_1)a_2 + (k_3 + t_3 + \tfrac{1}{2}(2k_2 + k_1 - k_1^2)t_1)a_3 + \sum_{j \geq 4}(k_j + t_j + t_1 k_{j-1})a_j.$$

Multiplication in A must be associative.
Let
$$x = \sum_{j \in \mathbf{N}} k_j a_j,\ y = \sum_{j \in \mathbf{N}} t_j a_j,\ z = \sum_{j \in \mathbf{N}} s_j a_j,\ xy = \sum_{j \in \mathbf{N}} r_j a_j,$$
$$(xy)z = \sum_{j \in \mathbf{N}} m_j a_j,\ yz = \sum_{j \in \mathbf{N}} u_j a_j,\ x(yz) = \sum_{j \in \mathbf{N}} v_j a_j.$$

We have

$m_1 = r_1 + s_1 = k_1 + t_1 + s_1,$
$v_1 = k_1 + u_1 = k_1 + t_1 + s_1,$
$m_2 = r_2 + s_2 + s_1 r_1 = k_2 + t_2 + t_1 k_1 + s_2 + s_1(k_1 + t_1) = k_2 + t_2 + t_1 k_1 + s_2 + s_1 k_1 + s_1 t_1,$
$v_2 = k_2 + u_2 + u_1 k_1 = k_2 + t_2 + s_2 + s_1 t_1 + (t_1 + s_1)k_1 = k_2 + t_2 + s_2 + s_1 t_1 + t_1 k_1 + s_1 k_1,$
$m_3 = r_3 + s_3 + \tfrac{1}{2}(2r_2 + r_1 - r_1^2)s_1 = r_3 + s_3 + r_2 s_1 + \tfrac{1}{2} r_1 s_1 - \tfrac{1}{2} r_1^2 s_1 =$
$k_3 + t_3 + \tfrac{1}{2}(2k_2 + k_1 - k_1^2)t_1 + s_3 + (k_2 + t_2 + t_1 k_1)s_1 + \tfrac{1}{2}(k_1 + t_1)s_1 - \tfrac{1}{2}(k_1^2 + t_1^2 + 2k_1 t_1)s_1 =$
$k_3 + t_3 + k_2 t_1 + \tfrac{1}{2} k_1 t_1 - \tfrac{1}{2} k_1^2 t_1 + s_3 + k_2 s_1 + t_2 s_1 + t_1 k_1 s_1 + \tfrac{1}{2} k_1 s_1 + \tfrac{1}{2} t_1 s_1 - \tfrac{1}{2} k_1^2 s_1 - \tfrac{1}{2} t_1^2 s_1 -$
$k_1 t_1 s_1 = k_3 + t_3 + k_2 t_1 + \tfrac{1}{2} k_1 t_1 - \tfrac{1}{2} k_1^2 t_1 + s_3 + k_2 s_1 + t_2 s_1 + \tfrac{1}{2} k_1 s_1 + \tfrac{1}{2} t_1 s_1 - \tfrac{1}{2} k_1^2 s_1 - \tfrac{1}{2} t_1^2 s_1,$
$v_3 = k_3 + u_3 + \tfrac{1}{2}(2k_2 + k_1 - k_1^2)u_1 = k_3 + u_3 + k_2 u_1 + \tfrac{1}{2} k_1 u_1 - \tfrac{1}{2} k_1^2 u_1 =$
$k_3 + t_3 + s_3 + \tfrac{1}{2}(2t_2 + t_1 - t_1^2)s_1 + k_2 t_1 + k_2 s_1 + \tfrac{1}{2} k_1 t_1 + \tfrac{1}{2} k_1 s_1 - \tfrac{1}{2} k_1^2 t_1 - \tfrac{1}{2} k_1^2 s_1 =$
$k_3 + t_3 + s_3 + t_2 s_1 + \tfrac{1}{2} t_1 s_1 - \tfrac{1}{2} t_1^2 s_1 + k_2 t_1 + k_2 s_1 + \tfrac{1}{2} k_1 t_1 + \tfrac{1}{2} k_1 s_1 - \tfrac{1}{2} k_1^2 t_1 - \tfrac{1}{2} k_1^2 s_1,$
$m_4 = r_4 + s_4 + s_1 r_3 = k_4 + t_4 + t_1 k_3 + s_4 + s_1(k_3 + t_3 + \tfrac{1}{2}(2k_2 + k_1 - k_1^2)t_1)) =$
$k_4 + t_4 + t_1 k_3 + s_4 + s_1 k_3 + s_1 t_3 + k_2 s_1 t_1 + \tfrac{1}{2} k_1 s_1 t_1 - \tfrac{1}{2} k_1^2 s_1 t_1,$
$v_4 = k_4 + u_4 + u_1 k_3 = k_4 + t_4 + s_4 + s_1 t_3 + (t_1 + s_1)k_3 = k_4 + t_4 + s_4 + s_1 t_3 + t_1 k_3 + s_1 k_3.$

And in general,

$m_{j+1} = r_{j+1} + s_{j+1} + s_1 r_j = k_{j+1} + t_{j+1} + s_{j+1} + t_1 k_j + s_1(k_j + t_j + t_1 k_{j-1}) =$
$k_{j+1} + t_{j+1} + s_{j+1} + t_1 k_j + s_1 k_j + s_1 t_j + s_1 t_1 k_{j-1},$
$v_{j+1} = k_{j+1} + u_{j+1} + u_1 k_j = k_{j+1} + t_{j+1} + s_{j+1} + s_1 t_j + (t_1 + s_1)k_j = k_{j+1} + t_{j+1} + s_{j+1} + s_1 t_j + t_1 k_j + s_1 k_j.$

Note that if $a_j = 0$, then $a_{j+1} = a_j * a = 0$, and hence $a_{j+k} = 0$ for all positive integer k.

Suppose that $2a_1 \neq 0$. Then we have:

$a_2 a_1 = 3a_1 + 2a_2 - a_3,\ ((2a_1)a_1)a_1 = (3a_1 + 2a_2 - a_3)a_1 = 4a_1 + 5a_2 - 2a_3 - a_4,$
$a_1 a_1 = 2a_1 + a_2,\ (2a_1)(a_1 a_1) = (2a_1)(2a_1 + a_2) = 4a_1 + 5a_2 - 2a_3.$

Thus, $(2a_1)(a_1 a_1) - (2a_1)(a_1 a_1) = a_4$. Hence, if $a_4 \neq 0$, then the multiplication is not associative. It follows that, in this case, $a_4 = 0$. Then, also, $a_5 = a_4 * a = 0$, and, furthermore, $a_j = 0$ for all $j > 5$. We can see that, in this case, A is nilpotent in the sense of Smoktunowizc, and more precisely, $A \in \mathcal{N}_S(4, 3)$.

If we suppose that $ma_1 = 0$ for some positive integer v, then

$$ma_2 = m(a_1 * a_1) = a_1 * a(ma_1) = 0,\ ma_3 = m(a_2 * a_1) = a_2 * a(ma_1) = 0,$$

and $ma_k = 0$ for all of positive integer k.

Now, consider the case when $2a_1 = 0$. Then $2a_2 = 0$ and $2a_k = 0$ for all of positive integer k. Thus, in this case, the additive group of A is elementary abelian 2 – group. We have:

$$a_2a_1 = a_1 + a_2 + a_3, (a_2a_1)a_1 = (a_1 + a_2 + a_3)a_1 = 2a_1 + 2a_2 + a_3 + a_4 = a_3 + a_4,$$
$$a_1a_1 = 2a_1 + a_2 = a_2, (a_2)(a_1a_1) = a_2 a_2 = 2a_2 = 0.$$

Equality $(a_2a_1)a_1 = (a_2)(a_1a_1)$ implies that $a_3 + a_4 = 0$ or $a_3 = a_4$. Then

$$a_3 = a_4 = a_3 * a_1 = a_3a_1 - a_3 - a_1 = a_3.$$

It follows that

$$a_3a_1 = a_3 + a_1 + a_3 = 2a_3 + a_1 = a_1.$$

The fact that by multiplication A is a group implies that $a_3 = 1 = 0$. Then, also, $a_j = 0$ for all $j > 3$. By **Proposition 2.1,** every element of A has a form $k_1a_1 + k_2a_2$, where $k_1, k_2$ are integers. Moreover, if $x = k_1a_1 + k_2a_2$, $y = t_1a_1 + t_2a_2$ are elements of A, $k_1, k_2, t_1, t_2 \in \mathbf{Z}$, then

$$x * y = t_1k_1a_2 \text{ and } xy = (k_1 + t_1)a_1 + (k_2 + t_2 + t_1k_1)a_2.$$

Furthermore, A is nilpotent in the sense of Smoktunowizc. More precisely, $A \in \mathcal{N}_S(3, 3)$.

## PROOF OF THEOREM B2.

Clearly, by addition, D is an abelian group. Furthermore, let

$$x = (k_1, k_2, k_3), y = (t_1, t_2, t_3), z = (s_1, s_2, s_3), xy = (r_1, r_2, r_3),$$
$$(xy)z = (m_1, m_2, m_3), yz = (u_1, u_2, u_3), x(yz) = (v_1, v_2, v_3).$$

We have:

$m_1 = r_1 + s_1 = k_1 + t_1 + s_1,$
$v_1 = k_1 + u_1 = k_1 + t_1 + s_1,$
$m_2 = r_2 + s_2 + s_1r_1 = k_2 + t_2 + t_1k_1 + s_2 + s_1(k_1 + t_1) = k_2 + t_2 + t_1k_1 + s_2 + s_1k_1 + s_1t_1,$
$v_2 = k_2 + u_2 + u_1k_1 = k_2 + t_2 + s_2 + s_1t_1 + (t_1 + s_1)k_1 = k_2 + t_2 + s_2 + s_1t_1 + t_1k_1 + s_1k_1,$
$m_3 = r_3 + s_3 + ½(2r_2 + r_1 - r_1^2)s_1 = r_3 + s_3 + r_2s_1 + ½ r_1s_1 - ½r_1^2s_1 =$
$k_3 + t_3 + ½(2k_2 + k_1 - k_1^2)t_1 + s_3 + (k_2 + t_2 + t_1k_1)s_1 + ½ (k_1 + t_1)s_1 - ½ (k_1^2 + t_1^2 + 2k_1t_1)s_1 =$
$k_3 + t_3 + k_2t_1 + ½k_1t_1 - ½k_1^2t_1 + s_3 + k_2s_1 + t_2s_1 + t_1k_1s_1 + ½k_1s_1 + ½t_1s_1 - ½k_1^2s_1 - ½t_1^2s_1 -$
$k_1t_1s_1 = k_3 + t_3 + k_2t_1 + ½k_1t_1 - ½k_1^2t_1 + s_3 + k_2s_1 + t_2s_1 + ½k_1s_1 + ½t_1s_1 - ½k_1^2s_1 - ½t_1^2s_1,$
$v_3 = k_3 + u_3 + ½(2k_2 + k_1 - k_1^2)u_1 = k_3 + u_3 + k_2u_1 + ½k_1u_1 - ½ k_1^2u_1 =$
$k_3 + t_3 + s_3 + ½(2t_2 + t_1 - t_1^2)s_1 + k_2t_1 + k_2s_1 + ½k_1t_1 + ½k_1s_1 - ½k_1^2t_1 - ½k_1^2s_1 =$
$k_3 + t_3 + s_3 + t_2s_1 + ½t_1s_1 - ½ t_1^2s_1 + k_2t_1 + k_2s_1 + ½k_1t_1 + ½k_1s_1 - ½k_1^2t_1 - ½k_1^2s_1.$

Thus, we obtain $m_1 = v_1$, $m_2 = v_2$, $m_3 = v_3$, which proves an equality $(xy)z = x(yz)$.

Thus, we can see that the multiplication is associative.
The identity element is $(0, 0, 0)$. Indeed,

$$(k_1, k_2, k_3)(0, 0, 0) = (k_1 + 0, 0 + k_2 + 0, k_3 + 0 + 0 ) = (k_1, k_2, k_3).$$

If $x = (k_1, k_2a, k_3)$, then $x^{-1} = (-k_1, k_1^2 - k_2, k_1k_2 + ½ (k_1^2 - k_1^3) - k_3)$.

Finally,

$x(y + z) = (k_1, k_2, k_3) ( (t_1, t_2, t_3) + (s_1, s_2, s_3) ) =$ ,
$x(y + z) = (k_1, k_2, k_3) ( t_1 + s_1, t_2 + s_2, t_3 + s_3 ) =$ ,
$( k_1 + t_1 + s_1, k_1(t_1 + s_1) + k_2 + t_2 + s_2, k_3 + t_3 + s_3 + ½ (2k_2 + k_1 - k_1^2)( t_1 + s_1) ) =$
$(k_1 + t_1 + s_1, k_1 t_1 + k_1 s_1 + k_2 + t_2 + s_2, k_3 + t_3 + s_3 + ½(2k_2 + k_1 - k_1^2)t_1 + ½(2k_2 + k_1 - k_1^2)s_1)$,

$xy + xz - x = (k_1, k_2, k_3) (t_1, t_2, t_3) + (k_1, k_2, k_3)(s_1, s_2, s_3) - (k_1, k_2, k_3) =$
$(k_1 + t_1, k_1 t_1 + k_2 + t_2, k_3 + t_3 + ½ (2k_2 + k_1 - k_1^2)t_1 ) +$
$(k_1 + s_1, k_1 s_1 + k_2 + s_2, k_3 + s_3 + ½ (2k_2 + k_1 - k_1^2)s_1 ) - (k_1, k_2, k_3) =$
$(k_1 + t_1, k_1 t_1 + k_2 + t_2, k_3 + t_3 + ½ (2k_2 + k_1 - k_1^2)t_1 ) +$
$(s_1, k_1 s_1 + s_2, s_3 + ½ (2k_2 + k_1 - k_1^2)s_1 ) - (k_1, k_2, k_3) =$
$(k_1 + t_1 + s_1, k_1 t_1 + k_1 s_1 + k_2 + t_2 + s_2, k_3 + t_3 + s_3 + ½(2k_2 + k_1 - k_1^2)t_1 + ½(2k_2 + k_1 - k_1^2)s_1)$,

so that $x(y + z) = xy + xz - x$.

This shows that D is a left brace.

Let A be an arbitrary left brace such that $A^3 = < 0 > = A^{(4)}$. Put $a_1 = a$, $a_2 = a * a = a_1 * a$, $a_3 = a_2 * a$. Since $A^{(3)} \neq < 0 >$, then $a_3 \neq 0$. By **Theorem A,** every element of A has a form $k_1 a_1 + k_2 a_2 + k_3 a_3$, where $k_1, k_2, k_3$ are integers.

Define the mapping f: D $\longrightarrow$ A by the rule $f(k_1, k_2, k_3) = k_1 a_1 + k_2 a_2 + k_3 a_3$. Using **Theorem A** again, we can show that f is an epimorphism.

The proof of Theorem B1 is similar.

L.A. Kurdachenko
Department of Algebra
Dnipro National University
Gagarin Prospect 70
Dnipro 10 (Ukraine)
e-mail: lkurdachenko@i.ua

I.Ya. Subbotin
Department of Mathematics
National University
5245 Pacific Concourse Drive
Los Angeles, CA 90045 (USA)
e-mail: isubboti@nu.edu